\documentclass[12pt]{article}

\usepackage{times}

\usepackage{epsf,latexsym,amssymb}
\usepackage{graphicx}
\usepackage{amsmath}  

\newcommand{\FIGDIR}{.}
\newcommand{\pic}[2]{\includegraphics[scale=#1]{\FIGDIR/#2}}
\newcommand{\picc}[2]{\begin{center}\pic{#1}{#2}\end{center}}

\newcommand{\comment}[1]{}
\newcommand{\edo}{\end{document}}

\newcommand{\R}{{\mathbb R}}  
\newcommand{\X}{{\mathbb X}}  
\newcommand{\Y}{{\mathbb Y}}  
 
\newcommand{\Rp}{\R_+} 
\newcommand{\Rnp}{\R^n_+} 
\newcommand{\Rnn}{\R^{n\times n}}
\newcommand{\U}{{\mathbb U}}  
\newcommand{\kon}{k_{\mbox{\tiny on}}}
\newcommand{\kontilde}{k_{\mbox{\tiny on}}^*}
\newcommand{\koff}{k_{\mbox{\tiny off}}}
\newcommand{\lF}{\lambda _{\mbox{\tiny F}}}
\newcommand{\vF}{v_{\mbox{\tiny F}}}
\newcommand{\wF}{w_{\mbox{\tiny F}}}
\newcommand{\rate}{\rho }
\newcommand{\lDIP}{\rate_{\mbox{\tiny DIP}}}
\newcommand{\btwo}{d} 
\newcommand{\tildevij}{\widetilde v_{ij}}

\def\twodigits#1{\ifnum#1<10 0\fi\the#1}
\makeatletter
  \newcommand\short[1]{\expandafter\\gamma obbletwo\number\numexpr#1\relax}
\makeatother
\usepackage{datetime2}

\newcommand{\abs}[1]{\left\vert #1 \right\vert}

\newtheorem{theorem}{Theorem}
\newtheorem{itlemma}{Lemma}
\newtheorem{itproposition}{Proposition}
\newtheorem{itcorollary}{Corollary}
\newtheorem{itremark}{Remark}
\newtheorem{itdefinition}{Definition}
\newtheorem{itexample}[itlemma]{Example}

\newenvironment{lemma}{\begin{itlemma}\rm}{\end{itlemma}} 
\newenvironment{remark}{\begin{itremark}\rm}{\end{itremark}} 
\newenvironment{corollary}{\begin{itcorollary}\rm}{\end{itcorollary}}
\newenvironment{proposition}{\begin{itproposition}\rm}{\end{itproposition}}
\newenvironment{definition}{\begin{itdefinition}\rm}{\end{itdefinition}}
\newenvironment{example}{\begin{itexample}\rm}{\end{itexample}}

\newcommand{\be}[1]{\begin{equation}\label{#1}}
\newcommand{\ee}{\end{equation}}
\newcommand{\bl}[1]{\begin{lemma}\label{#1}}
\newcommand{\ble}[1]{\begin{lemmaex}\label{#1}}
\newcommand{\br}[1]{\begin{remark}\label{#1}}
\newcommand{\bt}[1]{\begin{theorem}\label{#1}}
\newcommand{\bd}[1]{\begin{definition}\label{#1}}
\newcommand{\bp}[1]{\begin{proposition}\label{#1}}
\newcommand{\bc}[1]{\begin{corollary}\label{#1}}
\newcommand{\bex}[1]{\begin{example}\label{#1}}
\newcommand{\ec}{\mybox\end{corollary}}
\newcommand{\eex}{\mybox\end{example}}
\newcommand{\eem}{\mybox\end{example}}
\newcommand{\el}{\mybox\end{lemma}}
\newcommand{\er}{\mybox\end{remark}}
\newcommand{\et}{\qed\end{theorem}}
\newcommand{\ed}{\mybox\end{definition}}
\newcommand{\ep}{\mybox\end{proposition}}
\newcommand{\epr}{\end{proof}}
\newcommand{\bpr}{\begin{proof}}

\newcommand{\ecs}{\end{corollary}}
\newcommand{\eexs}{\end{example}}
\newcommand{\els}{\end{lemma}}
\newcommand{\ers}{\end{remark}}
\newcommand{\ets}{\end{theorem}}
\newcommand{\eds}{\end{definition}}
\newcommand{\eps}{\end{proposition}}
\newcommand{\halmos}{\rule{1ex}{1.4ex}}
\newcommand{\qed}{\hfill \halmos} 
\newcommand{\mybox}{\hfill $\Box$} 
\newcommand{\beq}{\begin{eqnarray}}
\newcommand{\eeq}{\end{eqnarray}}
\newcommand{\beqn}{\begin{eqnarray*}}
\newcommand{\eeqn}{\end{eqnarray*}}
\newcommand{\bi}{\begin{itemize}}
\newcommand{\ei}{\end{itemize}}
\newcommand{\ben}{\begin{enumerate}}
\newcommand{\een}{\end{enumerate}}

\newcommand{\bes}[1]{\begin{subequations}\label{#1}\begin{eqnarray}}
\newcommand{\ees}[1]{\end{eqnarray}\end{subequations}}

\newcommand{\uu}{u}
\newcommand{\vv}{v}
\newcommand{\ww}{w}
\newcommand{\Rate}{\overline{\rate}}
\newcommand{\fractext}[2]{#1/#2}

\newcommand{\minp}[2]{\begin{minipage}{#1\textwidth}#2\end{minipage}}

\newenvironment{proof}{\noindent {\em Proof}.\ }{\hspace*{\fill}$\halmos$\medskip}

\newcommand{\mypmatrix}[1]{\left(\begin{array}{cccccccccccc}#1\end{array}\right)}

\title{A concept of antifragility for dynamical systems}
\author{Eduardo D. Sontag\\
Northeastern University}

\begin{document}
\maketitle

\begin{abstract}
\footnotetext{timestamp:  {\small \today:\DTMcurrenttime}}

\noindent    
This paper defines antifragility for dynamical systems in terms of the convexity
of a newly introduced ``logarithmic rate''.  It
shows how to compute this rate for positive linear systems, and it
interprets antifragility in terms of pulsed alternations of extreme
strategies in comparison to average uniform strategies.

\end{abstract}

\section{Introduction}

Nassim Nicholas Taleb introduced in 2012 the concept of
``antifragility'' to refer to systems that \textit{benefit} from
uncertainty~\cite{taleb2012}. In contrast to \textit{robustness}
(resilience to uncertainty) being desirable, 
Taleb emphasized the advantages of devoting resources to ``placing
bets'' on unlikely events, provided that their payoff is large enough,
or more precisely that the expected return from a bet on extremes is
higher than the expected return from an intermediate bet, as will be
the case with convex payoff functions (Jensen's inequality, in
mathematical terms).
Much work has followed Taleb's original paper, including
formalizations in~\cite{taleb-douady2013} and recently a paper placing
these ideas in the context of mathematical
oncology~\cite{taleb-west2023}.

In this paper, we propose a definition of antifragility for dynamical
systems. The definition relies upon our introduction of a quantity, which
we call the ``logarithmic rate'' of an output, and which we write as
$\rate(u)$. The rate is a function of a parameter $u$ (and possibly
initial states).
Although antifragility is a very general idea that plays a role in areas
ranging from hedging strategies in finance to engineering~\cite{west-et-al2023},
to be concrete we will phrase our discussion in terms of infections
diseases or tumors. In that context, we may think of $u$ as quantifying a
dose of an antiviral or antibiotic to fight an infection, or a dose of
chemotherapy, immunotherapy, or targeted therapy in oncology, and we
may think of $\rate(u)$ as representing the rate of growth of an
infection or tumor, as discussed below. Antifragility will mean that
the function $\rate(u)$ is convex, if the objective is the
maximization of a reward (as seen by a tumor or a microorganism
carrying an infection) or that the function $\rate(u)$ is concave, if
the objective is the minimization of a cost (as seen by the infected
individual).
We will show how to compute $\rate(u)$ for the special, but important, case of
systems defined by positive linear dynamics, and we will focus on comparing
``pulsed'' versus ``uniform'' treatment protocols.

\subsection{Motivation}

\begin{figure}[ht]
  \minp{0.5}{\pic{0.4}{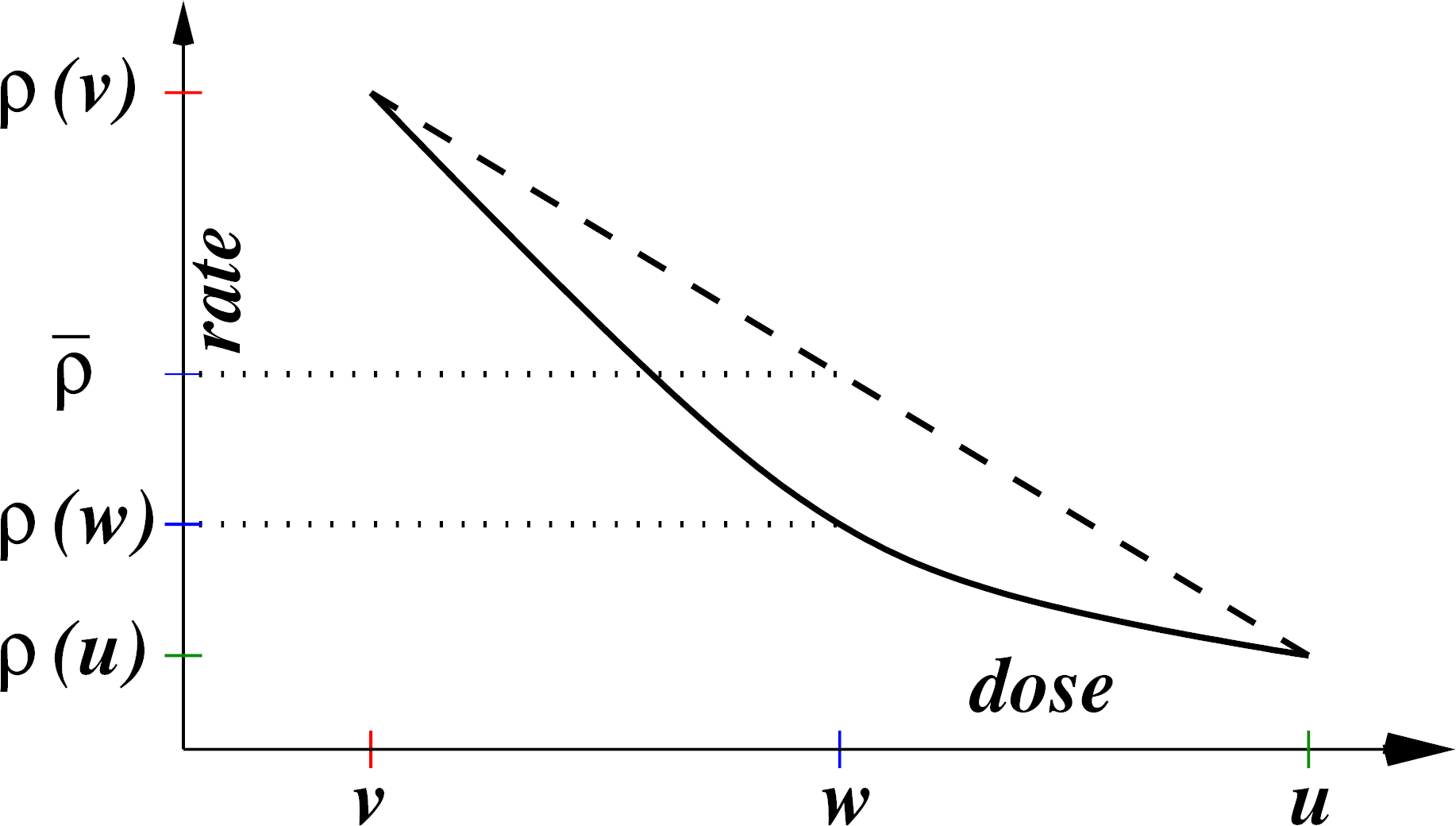}}%
  \minp{0.5}{\pic{0.4}{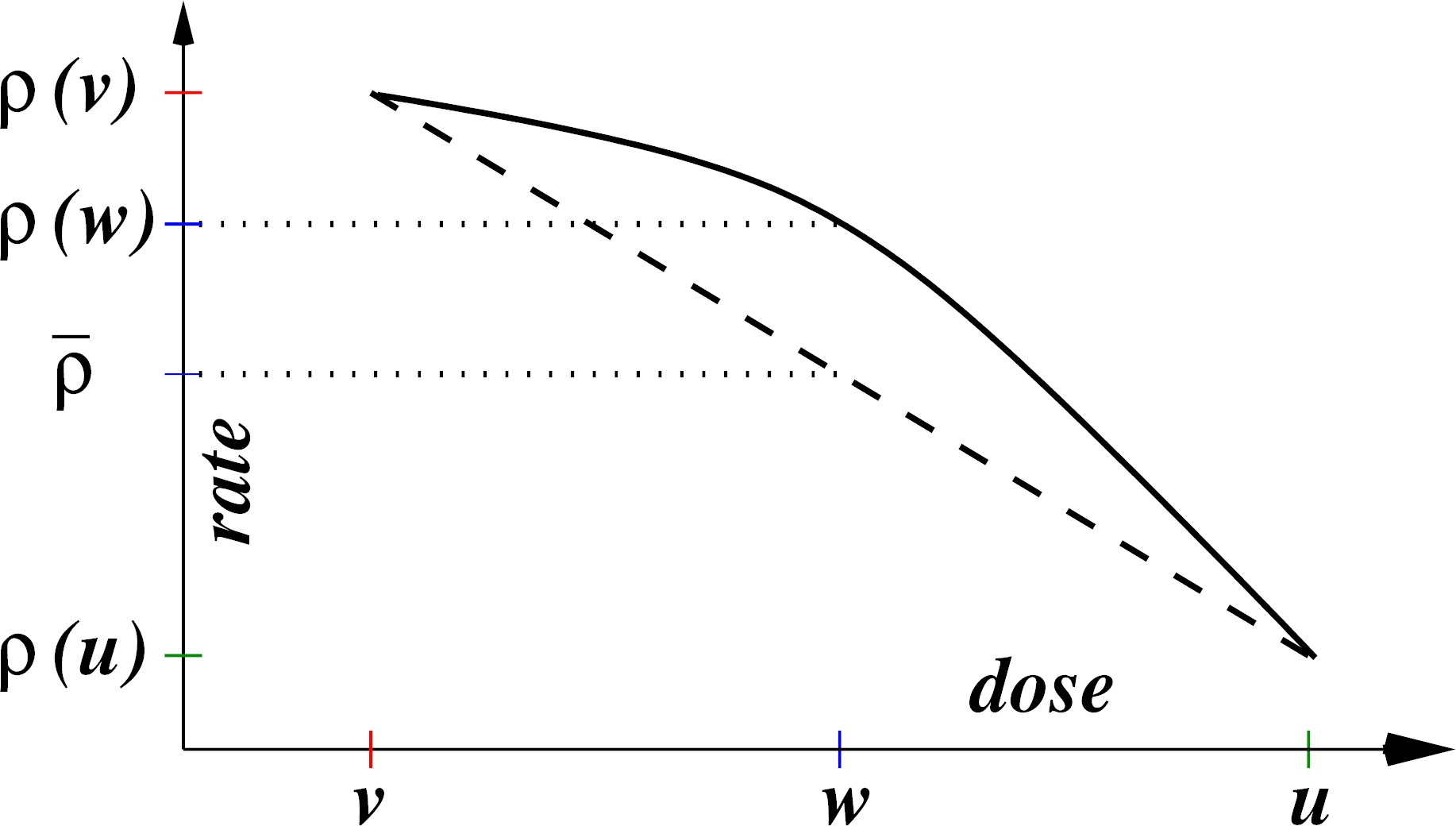}}
  \caption{Left: A convex function $\rate$, with dosages $\uu$, $\vv$, and
    $\ww = \frac{\uu+\vv}{2}$. Here $\rate(\ww)$ is smaller than $\Rate$.
    Right: Concave $\rate$; now $\rate(\ww)$ is larger than $\Rate$.}
  \label{fig:convex}
\end{figure}

To start the discussion, let us suppose that we have a 
cell population (consisting, for example, of cells infected by a
pathogen, or a type of tumor cells), and we wish to administer a
certain drug with the purpose of minimizing the population size.
Let us say that the size of the population at time $t$ is represented
by the scalar variable $x(t)$, which evolves according to the linear
differential equation 
\[
\dot x \;=\; \rate(u) \, x, \; x(0)=x_0\,.
\]
Here a dot indicates time derivative of $x=x(t)$ and $x_0$ is an
initial state, which we assume positive.
The function $\rate(u)$ quantifies the net growth rate when our drug is
given at a concentration quantified by $u$ (for example, $u=$ 10 in
units of mg per hour).
Note that $x(t) = e^{\rate(u)t}x_0$.
Our goal is that $x(t)$ should be small at some predetermined
future time $t=N$. 
When the drug is given at a high concentration, typically $\rate(u)$ will
be negative, so that the cell population tends to become extinct,
$x(t)=e^{\rate(u)N}x_0 \rightarrow  0$ as $N\rightarrow \infty $. For lower concentrations, on the
other hand, we might expect the rate of growth to be higher;
$\rate(u)$ could even be positive, in which case the population grows out of
control, $x(t)=e^{\rate(u)N}x_0\rightarrow  \infty $ as $N\rightarrow \infty $.
Clearly, we wish for $\rate(u)$ to be small. Assuming that a formula
for the function $\rate$ is known, such a minimization can be done
by solving a simple calculus problem. We could then administer the obtained
optimal dose $u$.
Fig~\ref{fig:convex} shows the graphs of two possible functions $\rho$.
If our purpose is that $\rate(u)$ should be as small as possible, we
can pick the rightmost points (both labeled $u$) in the graphs, that is
to say, the largest dose.

However, suppose now that, because of toxicity, or to minimize the
emergence of drug-induced resistance~\cite{greene_gevertz_sontag_resistance_2017},
we wish to consider more complicated therapies than simply
administering a constant drug concentration.
Is it better to use alternations of high and low doses or to use an
average dose? This question has been discussed in the context of
chemotherapy, where high or low variance treatment schedules may be
superior depending on the convexity of
$\rate(u)$~\cite{taleb-west2023,west-et-al2023}.  Specifically, let us
compare the following two scenarios for treatment.
 
\begin{figure}[ht]
  \picc{0.3}{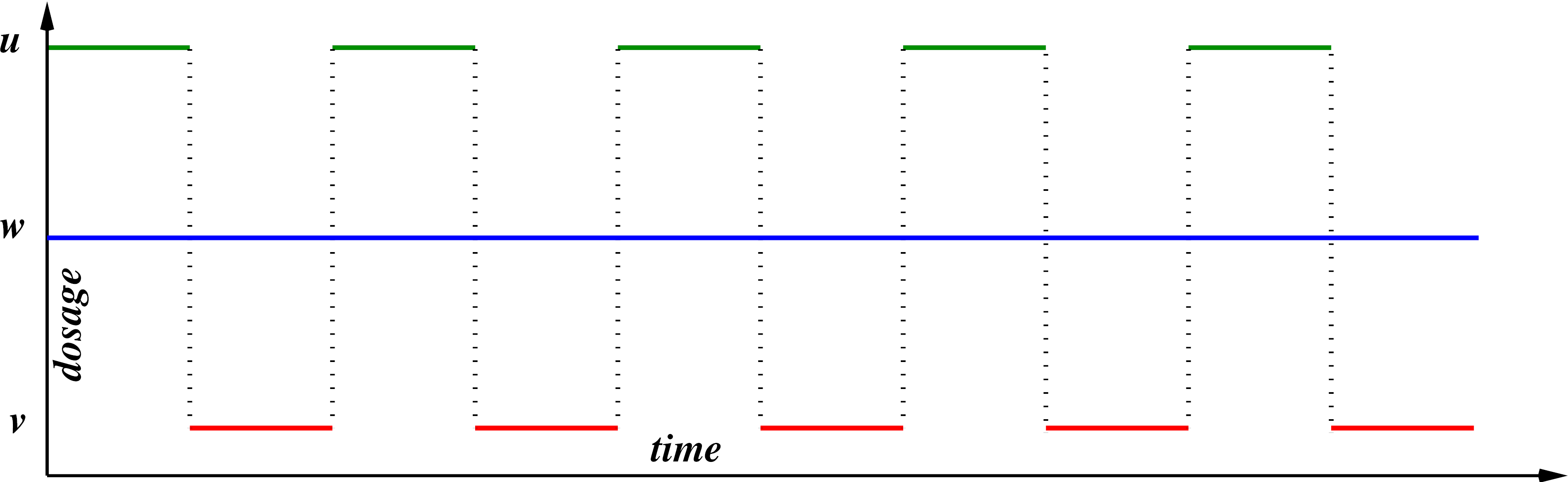}
  \caption{A pulsed protocol: iteration of high dose
    $\uu$ (green) followed by low dose $\vv$ (red). Shown also is a uniform protocol with average
    dose $\ww = \frac{\uu+\vv}{2}$ (blue).}
  \label{fig:treatment_schedule}
\end{figure}

In the first scenario, during a treatment period of total length $N$ we use
the following \textit{pulsed protocol}: we alternate between two drug dosages: a
higher dose $\uu$ followed by a lower dose $\vv$ (for example, $\vv=0$
would represent a ``drug holiday'' between treatments), each
applied for time $1/2$.
That is, dose $\uu$ is applied on the time intervals
\[
[0,0.5],
[1,1.5],
\ldots ,
[N-1,N-0.5]
\]
and dose $\vv$ is applied on the intervals
\[
[0.5,1],
[1.5,2],
\ldots 
[N-0.5,N]
\]
(see Fig~\ref{fig:treatment_schedule}).

\begin{figure}[ht]
  \picc{0.2}{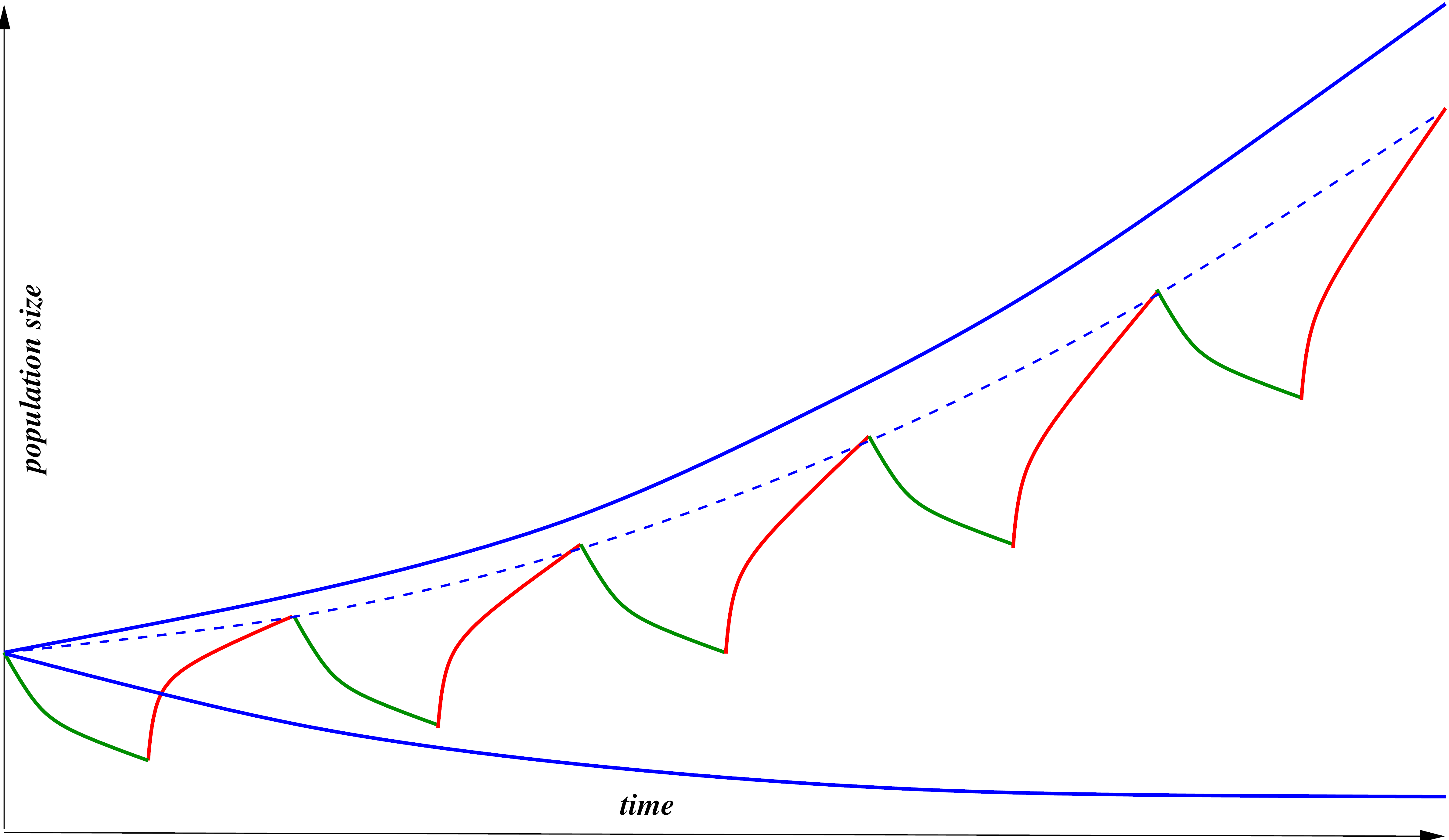}
  \caption{Population growth under pulsed protocol. Iteration of dose
    $\uu$ (green, with negative growth rate) followed by $\vv$ (red, positive
    growth rate). The resulting equivalent growth rate $\Rate$ is the average of
    the two rates; blue dashed line shows the equivalent growth under this
    rate, $e^{\Rate N}\, x_0$. A uniform protocol (solid blue lines)
    will provide a lower growth rate (perhaps negative) or a higher
    one, depending on the convexity of the rate function $\rate$.}
  \label{fig:protocols_result}
\end{figure}

As a result, in each period of length 1, the population size is first multiplied by
$e^{\rate(\uu)/2}$ and then by $e^{\rate(\vv)/2}$, so by the end, at
time, $t=N$, the  size of the population will be: 
\[
x_{\mbox{\tiny pulsed}}(N) \;=\; \left(
e^{\rate(\vv)/2}
e^{\rate(\uu)/2}
\right)^N
\;=\;
\left(
e^{\rate(\uu)/2+\rate(vv)/2}
\right)^N
\;=\; e^{\Rate N}\, x_0,
\]
where $\Rate$ is the average of the rates of $\uu$ and $\vv$,
$\Rate = \frac{1}{2}(\rate(\uu)+\rate(\vv))$.
See Fig~\ref{fig:protocols_result}.
Observe that the total amount of drug delivered during the even intervals is
$N\vv/2$ (there are $N$ intervals of length $1/2$, and in each of them the total
amount of drug is $\vv/2$).
Similarly, the total amount drug during the odd intervals is $N\uu/2$,
So the total amount of drug applied is
\[
\frac{N\uu}{2} + \frac{N\vv}{2} \;=\; N \frac{\uu+\vv}{2} \;=:\; N\ww\,.
\]

The second scenario is the following \textit{uniform protocol}: we
simply apply the constant dosage $\ww$ during the whole period $[0,N]$
(see Fig~\ref{fig:treatment_schedule}).
Now the size of the population at time $N$ is
\[
x_{\mbox{\tiny uniform}}(N) \;=\; e^{\rate(\ww)N}\,x_0 \,.
\]
The total amount of drug applied is $\ww N$, which is the same as in
the first scenario. However, if the rate function $\rate$ is not
linear, the outcome will generally be very different than
in the first scenario, 
A uniform protocol will result in a better (or at least not worse)
outcome --i.e.\ a lower $x(N)$-- than a pulsed protocol, if
$\rate(\ww)\leq \Rate$:
\be{eq:convex}
\rate\left(\frac{\uu+\vv}{2}\right) \;\leq \; \frac{1}{2} \bigg(\rate(\uu)+\rate(\vv)\bigg) \,,
\ee
which happens if the rate function $\rate$ is convex
(see Fig~\ref{fig:convex}, left).
Conversely, if the reverse inequality holds, which happens if the rate
function is concave, then a pulsed protocol is superior in decreasing
the population size
(see Fig~\ref{fig:convex}, right).

So far, the discussion has been limited to a very special situation, that of
a one-dimensional system, that is to say, there is only one state
variable describing the dynamics of the system. In a more complex
biological system, the effect of drugs can be complicated by phenomena
such as the emergence of therapy resistance, in which different
subpopulations, with different growth and death rates, react
differently to the drug. 
The main purpose of this paper is to provide a more complex setting
in which one can make a similar convexity-based comparison.
In this, we were motivated by the work
in~\cite{taleb-west2023,west-et-al2023} and closely related work in~\cite{dip}.

Still in the one-dimensional case, note that value of the initial condition is
not critical to the comparison of the sizes of the variable $x(t)$ for different
strategies, at least for large enough time. This is because the ratio
$e^{\rate_1t}x_1/e^{\rate_2t}x_2$ equals $ke^{(\rate_1-\rate_2)t}$,
for the constant $k=x_1/x_2$, and this ratio will be large if
$\rate_1>\rate_2$, and will tend to zero if the opposite inequality holds.
Thus we define a notion which is independent of the initial state.
(In section~\ref{sec:dip}, we briefly discuss what is ``large enough''
time.)
Since $x(t) = e^{\rate(u)t}x_0$, we can write $\ln x(t) =
\rate(u)t+\ln x_0$, which means that we may estimate the rate as
$\rate(u) = \fractext{\ln x(t)}{t} + \fractext{\ln x_0}{t}$. For large enough
times $t$, the effect of the initial condition is fades:
\[
\rate(u) \;=\; \lim_{t\rightarrow +\infty } \frac{\ln x(t)}{t}\,.
\]
We will provide an entropy-like quantity similar to $\rate$ for
arbitrary nonlinear dynamics, and show how to compute this quantity
for a class of linear systems $\dot x(t)=A(u)x(t)$, where now $x(t)$ is
vector valued and $A(u)$ is a matrix instead of the scalar $\rate(u)$.
When attempting to generalize to matrices the scalar case
discussed earlier, a key difficulty is that the formula
$e^{A(\vv)/2}e^{A(\uu)/2}=e^{A(\uu)/2+A(\vv)/2}$
is false
(though true in the special case that the matrices commute, $A(\uu)A(\vv)=A(\vv)A(\uu)$).
However, for the special class of ``positive'' systems, which arguably includes
all examples of biological interest, we are able to
derive an analogue of~\eqref{eq:convex} stated in terms of $\rate$,
thus allowing a comparison between
uniform and pulsed protocols using convexity of $\rho$ for
linear positive systems of dimension greater than one.
In the convex case, a pulsed protocol consisting
of sequential alternations of two drug concentrations results in a
larger growth rate than a uniform treatment using an average concentration.
Thus from the point of view of the tumor or infection, the pulsed protocol
results in a higher ``win'' compared to the uniform one.
From the patient's point of view, who wants to minimize instead of
maximize $\rho$, a constant strategy is best.
When the concavity is reversed, the opposite choices hold.

\subsection{Mathematical preliminaries}

Consider a parameterized system of differential equations as follows:
\be{eq:nonlinear}
\dot x = f(x,u), \;x(0) = x_0, \; y = h(x)\,.
\ee
Dot indicates time derivative, and arguments $t$ are not displayed.
Here $x=(x_1,\ldots ,x_n)$ is a state vector whose coordinates specify
population sizes, expression levels or differential expression levels
of molecular components such as mRNAs, or other quantities of
interest.
We assume that states $x(t)$ are required to belong to a prespecified
subset $\X$ of Euclidean space $\R^n$. this set $\X$ might include
constraints such as positivity of coordinates or maximal allowable
values for these coordinates. Here, $x_0\in \X$ is an initial condition.
The solution $x(t)=x(t,u,x_0)$ of~\eqref{eq:nonlinear} is assumed to
be unique and defined for all $t\geq 0$, taking values $x(t)\in \X$ for all $t$,
and $h:\X\rightarrow \Y$ is an output or readout map.
See~\cite{mct} for control-theory formalism and conditions
guaranteeing existence and uniqueness.
The symbol $u$ denotes a parameter that takes values in a convex
subset $\U$ of some Euclidean space.
We will denote by
\[
y(t,u,x_0)
\]
the output trajectory with input and initial state $x_0$. That is,
$y(t,u,x_0) = h(x(t,u,x_0))$, for all $t\geq 0$.

There are many applications of this setup. One cell biology
interpretation is as follows.  Each coordinate $x_i(t)$ of the state
vector quantifies the number (or volume) of cells of a certain type
``$i$'' at time $t$, $u$ represents a constant-in-time drug dosage
being applied during an interval of time, and $y(t)=h(x(t))$ is some
observable quantity, such as for example the total number of cells at
time $t$, $y(t)=\sum_{i=1}^nx_i(t)$.

\bd{def:rate}
The \textit{logarithmic rate} of the system is defined as follows:
\[
\rate\left(u,x_0\right) \,:=\; \limsup_{t\rightarrow +\infty } \frac{1}{t} \ln \abs{y(t,u,x_0)} 
\]
(when $y(t,u,x_0)=0$, we interpret the log as $-\infty $).
\ed

Note that when the limit exists, and if $y$ is positive, as will be the case
with the class of positive linear systems studied below, this is
simply
\[
\rate\left(u,x_0\right) \,:=\; \lim_{t\rightarrow +\infty } \frac{1}{t} \ln
y(t,u,x_0) \,. 
\]
In a typical application, we can think of $\rate$ is the slope of a
semi-log plot of cell number (or volume) vs.\ time.

The main motivation for this definition is as follows. Suppose that we
have a scalar ($n=1$) linear system $\dot x=\lambda x$ with measurement $y=cx$,
$x_0>0$, and $c>0$. Then $y(t) = c e^{\lambda t}x_0$, 
and therefore (since $cx_0\not= 0$ and thus its logarithm is well defined):
\[
\rate\left(u,x_0\right) \,=\; \lim_{t\rightarrow +\infty } \frac{1}{t} \ln y(t,u,x_0)
         \;=\; \lambda  + \lim_{t\rightarrow +\infty } \frac{\ln (cx_0)}{t} \;=\; \lambda 
\]
is independent of the initial state $x_0$.         
This is the rate of growth (or decay, if $\lambda <0$) of the output $y(t)$. 
When $\rate\left(u,x_0\right)$ is independent of $x_0$, as here, we
write simply $\rate(u)$.
Our main purpose here is to generalize this one-dimensional motivating
example to linear systems of higher dimension than one.

\br{rem:lyap_exponents}
A formula similar to that for $\rate$ is used when defining Lyapunov exponents
in chaos theory. A difference is that for Lyapunov exponents one looks at
differences between responses for two different initial states, which
would become an ``incremental'' variant of $\rate$, related to ideas
of logarithmic norms and contractive systems~\cite{aminzare_sontag_synchronization2014}.
\er

\subsection{Antifragility: reward maximization or cost minimization}

Intuitively, we want to define a system as ``antifragile'' if mixed
strategies are better than constant ones. For simplicity of
exposition, we assume from now on that $\rate$ does not depend on the
initial state. Mathematically:

\bd{def:antifragility}
The parameterized system~\eqref{eq:nonlinear} is
\textit{antifragile for reward maximization} if for
any $u_1,u_2\in \U$ and any nonnegative $\alpha _1,\alpha _2$ with $\alpha _1+\alpha _2=1$,
\[
\rate\left(\alpha _1u_1+\alpha _2u_2\right) \;\leq \;
         \alpha _1 \rate\left(u_1\right) + \alpha _2\rate\left(u_2\right),
\]
that is to say, the rate function $\rate$ is convex.
The parameterized system~\eqref{eq:nonlinear} is
\textit{antifragile for cost minimization} if the rate function
$\rate$ is concave (that is, the opposite inequality holds).
\ed

\br{rem:stochastic}
We motivated antifragility through sequential therapies.
A more standard probabilistic interpretation of antifragility is as follows.
Suppose that $\rate$ is a convex function, as the one shown in the
left panel of Fig~\ref{fig:convex}.  Now we think of $\rate(u)$ as the
``payoff'' of using an input $u$, thought of as a move in a game, an
investment strategy, or a bet on the outcome. Then, we get a better
expected payoff $\Rate$ when placing half of the time our bets on $u$
and half of the time on $v$, compared to placing all our bets at the
average $\ww=(\uu+\vv)/2$, since $\rate(\ww)\leq \Rate$.  If instead the
rate is concave, the expected payoff is better when using the average
$\ww$.
\er
We next discuss an important special case.

\section{Positive linear systems}

From now on, the set of states will be $\X =\Rnp$, the set of
$n$-tuples $(x_1,\ldots ,x_n)$ with positive coordinates.
In particular, the initial condition $x_0$ will have positive entries.
Outputs will take values in $\Y = \Rp$, the set of positive real
numbers.
These are of course reasonable assumptions when we are dealing with
biological populations, volumes, levels of protein expression, and so forth.
The main example of systems~\eqref{eq:nonlinear}
to be discussed here is as follows. We consider
\textit{positive linear systems}, meaning systems of the form
\be{eq:positive}
\dot x = A(u)x,\; x(0) = x_0, \; y = c x
\ee
where $A(u)\in \Rnn$ is a parameterized \textit{Metzler matrix}, meaning
\cite{farina2000} a matrix whose off-diagonal elements are nonnegative
(diagonal elements are allowed to be positive, zero, or negative; also
called an ``essentially nonnegative'' matrix in~\cite{horn-jackson}).
and $c$ is a row vector $(c_1,\ldots ,c_n)$ with positive entries.
It is known that for such systems, if the initial condition has
positive entries then the solution $x(t)$ will also have positive
entries for every $t>0$.
We will assume the following irreducibility condition: for some
$r\in \R$, and some positive integer $k$, the matrix $(rI+A)^k$ has
all entries positive. A sufficient condition for irreducibility to
hold is if the off-diagonal elements of $A$ are all strictly positive:
indeed, just add a sufficiently large $r$ to make all diagonal elements
positive; then the condition holds with $k=1$.
(Metzler irreducible matrices are closely related to a more general class
called primitive matrices, for which the properties stated below
are also true, but we do not need the more general concept.)

Using the Jordan canonical form, it is known that the general solution
of a system $\dot x=Ax$ is a sum of terms of the form 
\be{eq:gensol}
\beta _{ij} e^{\lambda _it} \left(\frac{t^{j-1}}{(j-1)!} v_{i,1} +
\frac{t^{j-2}}{(j-2)!} v_{i,2} + \ldots 
+ \frac{t}{1!} v_{i,j-1} + v_{i,j}\right)
\;=\;
\beta _{ij} e^{\lambda _it}\tildevij(t)
\ee
where the $\lambda _i$'s are eigenvalues of $A$ and the $v_{i,j}$'s are
eigenvectors or generalized eigenvectors of $A$.
Each $\tildevij(t)$ is a vector of polynomials in $t$.
The constants $\beta _{ij}$ (which may be complex) are obtained
from the initial condition $x_0$, which we assume to be positive, by
expanding $x_0$ in terms of eigenvectors and generalized eigenvectors of $A$.
The coefficients of this expansion are obtained from a linear
transformation of the coordinates of $x$ in the canonical basis, so
there are row vectors $q_{ij}$ such that $\beta _{ij}=q_{ij}x_0$.

We now specialize this general solution to irreducible Metzler
matrices $A$.
For such matrices, it is known~\cite{farina2000} that there is an
eigenvalue $\lF$ of $A$, called the Frobenius eigenvalue, which is
real and which dominates all other eigenvalues, meaning that all other
eigenvalues have real part less than $\lF$. Moreover, there is a
(real) eigenvector $\vF$ associated to $\lF$ which has all its entries
positive, and all other eigenvectors associated to $\lF$ are multiples
of $\vF$ ($\lF$ has algebraic multiplicity one). 
Therefore, among the expressions in~\eqref{eq:gensol} there is one of
the form $\beta  e^{\lF t}\vF$
with $\beta $ real, and all others have the form
$[q_{ij}x_0] e^{\lambda _it}\tildevij(t)$ with
$\Re \lambda _{i} < \lF$. Each such term can be written as 
$e^{\lF t} [q_{ij}x_0] e^{\mu _it}\tildevij(t)$ with $\Re \mu _i<0$.
Collecting all these other terms, we can write them together as
$e^{\lF t} W(t)x_0$, where $W(t)$ is an $n\times n$ matrix whose entries
are linear combinations of terms of the form $t^k e^{\mu _it}$ with $\Re
\mu _i<0$, and hence $W(t)\rightarrow 0$ as $t\rightarrow +\infty $.
Thus the general solution has the form
$x(t) = e^{\lF t} \left( \beta \vF + W(t)x_0 \right)$.
Moreover,
$\beta \vF$ is the projection of the initial
condition onto the eigenspace corresponding to $\lF$, $Px_0 = \beta  \vF$.
The matrix $P$ is the \textit{Perron projection}, which has the formula
$P = \vF \wF$, where $\wF$ is a row vector, a positive
left eigenvector of $A$ with eigenvalue $\lF$, i.e.\ $\wF A = \lF\wF$, 
and it is picked so that $\wF\vF=1$. In summary,
\be{eq:gensolpos}
x(t) \;=\; e^{\lF t} \left( P + W(t) \right) x_0\,.
\ee
Recall that our system is positive, so $x(t)$ has positive entries as
long as the initial condition was also positive.
Suppose now that $y(t)=y(t,u,x_0)=c\,x(t,u,x_0)$,
where $c=(c_1,\ldots ,c_n)$ with all $c_i>0$.
It follows that
\be{eq:output}
y(t) \;=\; e^{\lF t} (\kappa  + \mu (t))
\ee
where $\kappa  = c\,Px_0$ is positive since $P$, $c$, and $x_0$ all have
positive entries, and $\mu (t)=c\,W(t)x_0\rightarrow 0$ as $t\rightarrow \infty $.
Moreover, $y(t)$ is positive for every
$t>0$, so that the natural log $\ln y(t)$ is well defined.
Thus
\[
\lim_{t\rightarrow \infty } \frac{1}{t} \ln y(t)
\;=\; \lF + \lim_{t\rightarrow \infty } \frac{1}{t} \ln(\kappa +\mu (t)) \;=\; \lF\,.
\]
Let is write now $\lF(u)$ instead of $\lF$ in order to emphasize that
$A=A(u)$.
We have proved:

\bt{theo:linear-positive}
For irreducible positive systems,
$
\rate(u) \;=\; \lF(u)\,.
$
\et

This allows is to calculate $\rate(u)$ explicitly, as well as to
compute the effect of sequential strategies, as done below.

\section{Sequential inputs}

Suppose now that we consider two inputs, $u_1$ and $u_2$, and these
are consecutively applied: first input $u_1$ is applied for time
duration $\alpha t$, and then input $u_2$ is applied, for time duration
$\beta t$, where $\alpha $ and $\beta =1-\alpha $ are positive numbers.
So the total time is $t$.

Let us write $x(t,u_1,u_2,x_0)$ for the resulting state at time $t$,
and $y(t,u_1,u_2,x_0)=h(x(t,u_1,u_2,x_0))$.
Equivalently, we first solve the differential equation with initial
state $x_0$ and input $u_1$, for time $\alpha t$, and then use the
resulting state $x(\alpha t)$ as a new initial state and then solve
the differential equation for time $\beta t$, but now with input $u_2$.
This means that
\[
x(t,u_1,u_2,x_0) \;=\; x(\beta t,u_2,x(\alpha t,u_1,x_0)).
\]
We define, if the limit exists,
\[
\rate(u_1,u_2,\alpha ) \,:=\; \lim_{t\rightarrow +\infty } \frac{1}{t} \ln y(t,u_1,u_2,x_0) \,.
\]
We can think of this quantity as a logarithmic rate for a sequential
application of inputs.
(A longer periodic alternation of $u_1$ and $u_2$ would be defined in a
similar way.)
From now on, we specialize to irreducible positive linear systems.

Using~\eqref{eq:gensolpos}, we know that
\[
x(t) \;=\; e^{\lF(u_2) \beta t} \left( P(u_2) + W_2(t) \right) x_1(\alpha t) 
\]
where $\lF(u_2)$ and $P(u_2)$ are the Perron eigenvalue and projection
matrix, respectively, corresponding to input $u_2$,
$W_2(t)\rightarrow 0$ as $t\rightarrow +\infty $,
and $x_1(\alpha t)=x(\alpha t,u_1,x_0)$.
Once again using~\eqref{eq:gensolpos},
\[
x_1(\alpha t) \;=\; e^{\lF(u_1) \alpha t} \left( P(u_1) + W_1(t) \right) x_0
\]
where $\lF(u_1)$ and $P(u_1)$ are the Perron eigenvalue and projection
matrix, respectively, corresponding to input $u_1$,
and $W_1(t)\rightarrow 0$ as $t\rightarrow +\infty $.
It follows that
\beqn
x(t) &=&
e^{[\alpha \lF(u_1) + \beta \lF(u_2)]t}
\left( P(u_2) + W_2(t) \right) \left( P(u_1) + W_1(t) \right) x_0\\
&=&
e^{[\alpha \lF(u_1) + \beta \lF(u_2)]t}
\left( P(u_2) P(u_1) + W(t) \right) x_0
\eeqn
where
\[
W(t)\;=\; P(u_2)W_1(t) + W_2(t)P(u_1) + W_1(t)W_2(t) \;\rightarrow \; 0 \;
\mbox{as} \; t\rightarrow  +\infty \,.
\]
Therefore
\[
y(t) \;=\; e^{(\alpha \lF(u_1) + \beta \lF(u_2))t} (\kappa +\mu (t))
\]
where $\kappa  = cP(u_2)P(u_1)x_0$ is positive (since both matrices
$P(u_1)$ and $P(u_2)$
and the vectors $c$ and $x_0$ have all entries positive) and
$\mu (t) = c W(t) x_0 \rightarrow 0$ as $t\rightarrow +\infty $.
Thus
\[
\lim_{t\rightarrow \infty } \frac{1}{t} \ln y(t)
\;=\;
\alpha \lF(u_1) + \beta \lF(u_2) + \lim_{t\rightarrow \infty } \frac{1}{t} \ln(\kappa +\mu (t))
\;=\;
\alpha \lF(u_1) + \beta \lF(u_2)\,.
\]
We have therefore proved:

\bt{th:alternation}
For irreducible positive systems, $\rate(u_1,u_2,\alpha ) = \alpha \lF(u_1) + (1-\alpha )\lF(u_2)$.
\et

The significance of this result is that, if we wish to compare the
``gain'' from an average input $u$ to that of using sequential inputs
$(u_1,u_2)$, the second will give a better (or equal) result, at least for large
total time $t$, if and only if
\[
\alpha \lF(u_1) + (1-\alpha )\lF(u_2) \;\geq \; \lF(u)\,,
\]
in other words, if and only if the rate function $\lF$ is convex.

\section{Flux-growth interpretation}
\label{sec:dip}

Observe that any $n\times n$ matrix $A$ with positive off-diagonals elements
can be written in the following form:
\[
\mypmatrix{
-\sum_{j\not= 1}a_{j1} + b_1 & a_{12}                  & \cdots & a_{1n}\\
                 a_{21} & -\sum_{j\not= 2}a_{j2} + b_2 & \cdots & a_{2n}\\
                 \vdots  & \vdots                 & \ddots & \vdots\\
                 a_{n1} & a_{n2} & \cdots & -\sum_{j\not= n}a_{jn} + b_n}
\]
with $a_{ij}>0$ for every $i\not= j$, and some $b_i$'s: to see this, simply define
$b_i:=a_{ii}+\sum_{j\not= i}a_{ij}$.

This trivial observation allows the following interpretation: the
matrix $A$ can be seen as describing fluxes $a_{ij}$ from
cells of type $j$ into cells of type $i$: the negative diagonal
elements combine the sum of all fluxes ``out of'' the cell $i$ with
the net growth rate $b_i$ (negative or positive) of cells of type $i$.

\br{rem:fast}
Suppose now that all off-diagonal elements are positive and
\textit{large}, corresponding to high fluxes between compartments.
In that case, the term $\mu (t)$ in~\eqref{eq:gensolpos} will tend to zero
very fast, and thus comparisons based on $y(t)$ instead of rates will
still be valid for ``small'' times.
Let us sketch the argument.
Introduce the matrix $A_0$ obtained from $A$ by setting all
$b_i$ to zero.  The sum of each column of $A_0$ is zero, implying
that the matrix $A_0$ has an eigenvalue $\lambda =0$ corresponding to the
right eigenvector $\mathbf{1}=[1,1\ldots ,1]^\intercal$, and all the other
eigenvalues $\lambda _2,\ldots ,\lambda _n$ of $A_0$ have negative real parts
(irreducible Metzler matrix). These real parts will be large in
magnitude, since we can think of $A_0$ as a rescaling by a large
scalar of a matrix with $n-1$ eigenvalues with negative real part and
one eigenvalue at zero.  Note that $A = A_0+B$, where $B$ is a
diagonal matrix which has the $b_i$'s on its diagonal.
Thus $A$ is a bounded perturbation of $A_0$ ($B$ is fixed, but the
other entries of $A$ can be though of as tending to infinity).
By continuity of eigenvalues on matrix entries (see for example
Corollary A.4.3 in~\cite{mct}), the eigenvalues of $A$ will be near
those of $A_0$. Therefore, the gap between $\lF$ and the real part of the
next larger eigenvalue remains large even after the perturbation. This means that the term $\mu (t)$
converges to zero very fast, as claimed.
In the next section, we make this argument explicit for the
two-dimensional case.
\er

\subsection{An example: two-dimensional systems}

In the special case $n=2$ of two cell types, this would be written as:
\be{eq:2dgeneral}
A \;=\;
\mypmatrix{
  -\kon + b_1 & \koff\\
   \kon       & -\koff + b_2}
\ee
where $\kon$ and $\koff$ are positive numbers and $b_1$ and $b_2$ are
arbitrary numbers. 
Thus $\kon$ is the rate at which cells of type 1 transition to cells
of type 2, $\koff$ is the rate at which cells of type 2 transition
to cells of type 1, and $b_1,b_2$ are the respective growth rates of
cells of types 1 and 2.

In general, all four numbers might depend on $u$.
An even more specific example 
would have the $u$ dependence in just the $\kon$ term, so we replace
$\kon$ by $\kon = \kontilde u$ in~\eqref{eq:2dgeneral} and assume that
the remaining entries do not depend on $u$.
For notational simplicity in future calculations, let us rewrite
such a matrix as follows: 
\be{eq:2d}
A(u) \;=\;
\mypmatrix{
  b -a u  & ak\\
   au     & d-ak}
\ee
where we are letting
$a = \kontilde$, $b=b_1$, $d=b_2$, and $k=\koff/\kontilde$.
This represents a system in which
cells of type 1 transition to cells of type 2 at a rate $\kontilde u$
that is proportional to the drug concentration $u$. Conversely, cells of type 2
transition back to type 1 at a rate that is independent of the drug
concentration. (For simplicity, we use the term ``drug concentration''
but this term typically represents the actual effect of the drug,
which could be a phamacokinetically derived Michaelis-Menten or Hill
function of the dose being applied to the system.) 

\subsection{Connection to drug-induced proliferation rate}

The paper~\cite{dip} introduced the idea of using the \textit{drug-induced proliferation rate} (``DIP'')
as a metric to quantify and compare drug effects, in the $n=2$ case.
We show next that the DIP rate is the same as the logarithmic rate
of the system, under the same approximation as used in~\cite{dip} to
derive it.

The paper defines the DIP as ``the steady-state rate of proliferation of a
cell population in the presence of a given concentration of drug.''
Specifically, the $A$ matrix is of the form shown in
equation~\eqref{eq:2d}, and $c = (1, 1)$.
The paper provides the approximate formula
\be{eq:dip}
y(t) \;\approx\; e^{\lDIP(u) t} x_0,
\ee
where
\[
\lDIP(u) \;=\; \frac{bk + du}{k+u} \,.
\]
The formula is derived by assuming a fast exchange between cell types,
giving an equilibrium
approximation
$\kontilde u x_1(t) = \koff x_2(t)$, and substituting into the two
differential equations (see Supplementary Note equations S14-S20 in~\cite{dip}).
This amounts to asking that both $\kontilde$ and $\koff$ ($a$ and $ak$) are
large while keeping the ratio $k = \kontilde/\kon$ constant.
In other words, the assumption is that $a\rightarrow +\infty $ while maintaining
$c,d,k$ constant.
We write $\rate(a,c,k,u)$ to emphasize the dependence of $\rate$ on
the parameters.

\bt{theo:equal-DIP}
$\lim_{a\rightarrow +\infty }\rate(a,c,k,u) = \lDIP(u)$.
\ets

\bpr
A little calculus exercise (see below) shows that for any positive
numbers $p,q$:
\be{eq:lhopital}
\lim_{x \to \infty} \left( -x + \sqrt{x^2 - px + q} \right)
\;=\;
-p/2\,.
\ee
Now note that
\[
\mbox{trace}\,(A) \;=\; T \;=\; b+{\btwo}-a(k+u)
\]
and
\[
\mbox{det}\,(A) \;=\; D \;=\; b{\btwo} -a(bk+{\btwo}u).
\]
Let us define:
\[
p \;:=\; -\frac{4\,\left(bk+{\btwo}u\right)}{k+u}
\]
\[
q \;:=\;\frac{4\left(kb^2+u{\btwo}^2\right)}{k+u}
\]
\[
x \;:=\; a\left(k+u\right) - (b+{\btwo}).
\]
A substitution and simplification shows that
\[
T^2 - 4D \;=\; x^2 - px +q \,.
\]
Therefore, the largest eigenvalue of $A$ is:
\[
\lambda(a,b,k,u) \;=\; (1/2)\left(T - \sqrt{T^2-4D}\right) \;=\; (1/2) \left(-x + \sqrt{x^2-px+q}\right)\,.
\]
Since $x\to \infty$ is equivalent to $a\to \infty$, and
applying~\eqref{eq:lhopital} to our $p$, we conclude that:
\[
 \lim_{x \to \infty} \lambda(a,b,k,u)  \;=\; 
\frac{bk+{\btwo}u}{k+u}
\]
as claimed.

For completeness, let us show how to establish~\eqref{eq:lhopital}. 
let $L:= \lim_{x \to \infty} \left( -x + \sqrt{x^2 - px + q} \right)$.
We will show that $L=-p/2$.
The given expression is of the
indeterminate form ``$-\infty + \infty$'' so we need to apply
L'H\^opital's rule.
We first multiply and divide by the conjugate:
\begin{align*}
L &= \lim_{x \to \infty} \frac{\left(-x + \sqrt{x^2 - px + q}\right) \left(-x - \sqrt{x^2 - px + q}\right)}{-x - \sqrt{x^2 - px + q}}\\
  &= \lim_{x \to \infty} \frac{px - q}{-x - \sqrt{x^2 - px + q}}\,.
\end{align*}
Now the expression is in the form $\frac{\infty}{\infty}$, so we apply
L'H\^opital's rule by differentiating the numerator and the denominator.
The derivative of the numerator $px - q$ is $p$.
The derivative of the denominator $-x - \sqrt{x^2 - px + q}$ is computed as:
\begin{align*}
  \frac{d}{dx} \left( -x - \sqrt{x^2 - px + q} \right)
&= -1 - \frac{1}{2\sqrt{x^2 - px + q}} \cdot (2x - p) \\
&= -1 - \frac{x - \frac{p}{2}}{\sqrt{x^2 - px + q}}\,.
\end{align*}
Thus, applying L'H\^opital's Rule gives:
\[
L = \lim_{x \to \infty} \frac{p}{-1 - \frac{x + \frac{p}{2}}{\sqrt{x^2 - px + q}}}\,.
\]
As $x \to \infty$, the term $\frac{x - \frac{p}{2}}{\sqrt{x^2 - px +
    q}}$ approaches $1$, so the denominator converges to -2.
Thus, the limit becomes:
\[
L = \frac{p}{-2} = -\frac{p}{2}
\]
as claimed.
\epr

\subsection*{Acknowledgments}

The author thanks Dr.\ Jeffrey West for very productive and informative
discussions, as well as making his papers and preprints available.
This research was supported in part by grants ONR N00014-21-1-2431
and AFOSR FA9550-21-1-0289.

\edo